\documentclass[12]{amsart}
\usepackage{amsfonts}
\newtheorem{thm}{Theorem}[section]
\newtheorem{cor}[thm]{Corollary}
\newtheorem{lem}[thm]{Lemma}

\newtheorem{prop}[thm]{Proposition}
\newtheorem{exa}[thm]{Example}

\theoremstyle{definition}
\newtheorem{defn}[thm]{Definition}
\usepackage{amscd}
\newtheorem{rem}[thm]{Remark}

\numberwithin{equation}{section}

\begin{document}
\title[On strongly $g(x)$-clean rings ]{On strongly $g(x)$-clean rings}
\author[Fan and Yang]{Lingling Fan $^{\dagger}$ and Xiande Yang$^{\ddagger}$}
\date{\today}\maketitle
\centerline{\tiny {$\dagger$Department of Mathematics and
Statistics, Memorial University of Newfoundland, St.John's, A1C 5S7,
Canada}}
 \centerline{\tiny {Email address: lingling@math.mun.ca}}

 \centerline{\tiny {$\ddagger$Department of Mathematics, Harbin Institute of Technology, Harbin, 150001, China}}
\centerline{\tiny {Email address: xiande.yang@gmail.com }}
\footnote{$^{\ddagger}$: Corresponding author}
\begin{abstract}
Let $R$ be an associative ring with identity, $C(R)$ denote the
center of $R$, and $g(x)$ be a polynomial in the polynomial ring
$C(R)[x]$. $R$ is called strongly $g(x)$-clean if every element $r
\in R$ can be written as $r=s+u$ with $g(s)=0$, $u$ a unit of $R$,
and $su=us$. The relation between strongly $g(x)$-clean rings and
strongly clean rings is determined, some general properties of
strongly $g(x)$-clean rings are given, and strongly $g(x)$-clean
rings generated by units are discussed.

\medskip


\end{abstract}

\bigskip
\baselineskip=18pt

\bigskip

\section{Introduction}
Let $R$ be an associative ring with the group of units $U(R)$. $R$
is called {\it clean} if for every element $r \in R $, $r=e+u$ with
$e^2=e \in R$ and $u\in U(R)$ \cite{N77} and $R$ is called {\it
strongly clean} if in addition, $eu=ue$ \cite{N99}.

Let $C(R)$ denote the center of a ring $R$ and $g(x)$ be a
polynomial in $C(R)[x]$. Camillo and Sim\'{o}n \cite{CS02} say $R$
is \emph {$g(x)$-clean} if for every element $r \in R$, $r=s+u$ with
$g(s)=0$ and $u \in U(R)$. If $V$ is a countable dimensional vector
space over a division ring $D$, Camillo and Sim\'{o}n proved that
$End(_{D}V)$ is $g(x)$-clean if $g(x)$ has two distinct roots in
$C(D)$ \cite{CS02}. Nicholson and Zhou generalized Camillo and
Sim\'{o}n's result by proving that $End(_{R}M)$ is $g(x)$-clean if
$_{R}M$ is a semisimple $R$-module and $g(x)\in (x-a)(x-b)C(R)[x]$
where $a, b\in C(R)$ and $b, b-a \in U(R)$ \cite{NZ05}. \cite{FY,WC}
completely determined  the relation  between clean rings and
$g(x)$-clean rings independently. What is the relation between
strongly clean rings and $g(x)$-clean rings?

In this paper, we continue this topic. In Section $2$, we define
strongly $g(x)$-clean rings and determine the relation between
strongly $g(x)$-clean rings and strongly clean rings; in Section
$3$, some general properties of strongly $g(x)$-clean rings are
given; and in Section $4$, some classes of strongly $g(x)$-clean
rings generated by units are discussed.

Throughout the paper,  ${\mathbb T}_n(R)$ denotes the upper
triangular matrix ring of order $n$ over $R$, ${\mathbb N}$ denotes
the set of all positive integers, and ${\mathbb Z}$ represents the
ring of integers.

\medskip \quad

\section{strongly $g(x)$-clean rings vs strongly clean rings}
\begin{defn} Let $g(x)\in C(R)[x]$ be a fixed polynomial. An element $r\in R$ is {\it strongly} $g(x)$-{\it clean}
if $r=s+u$ with $g(s)=0$, $u\in U(R)$, and $su=us$. $R$ is {\it
strongly} $g(x)$-{\it clean} if every element of $R$ is strongly
$g(x)$-clean. \end{defn}

Strongly clean rings are exactly strongly $(x^2-x)$-clean rings.
However, there are strongly $g(x)$-clean rings which are not
strongly clean and vice versa:

 Let ${\mathbb Z}_{(p)} =\{\frac{m}{n}\in {\mathbb Q}: {\rm gcd}(p,
 n)=1 \text{ and } p \text{ prime  }\}$ be the localization of ${\mathbb Z}$ at the prime ideal $p{\mathbb Z}$
 and $C_3$ be the cyclic group of order $3$.

\begin{exa}\label{exa:2.1} Let $R$ be a
commutative local or commutative semiperfect ring with $2 \in U(R)$.
By the proof of \cite[Theorem 2.7]{WC}, $RC_{3}$ is strongly
$(x^{6}-1)$-clean. In particular, ${\mathbb Z}_{(7)}C_3$ is a
strongly $(x^6-1)$-clean ring. Furthermore, by \cite[Example
2.2]{FY}, ${\mathbb Z}_{(7)}C_3$ is strongly $(x^4-x)$-clean.
However, ${\mathbb Z}_{(7)}C_3$ is not strongly clean \cite[Example
1]{HN01}.
\end{exa}

\begin{exa}\label{exa:3}
Let $R={\mathbb Z}_{(p)}$ and $g(x)=(x-a)(x^{2}+1) \in C(R)[x]$.
Then $R$ is strongly clean but by a easy verification we know $R$ is
not strongly $g(x)$-clean. Let $R$ be a boolean ring with more than
two elements with $c\not=0, 1$. Then $R$ is strongly clean but $R$
is not strongly $g(x)=(x+1)(x+c)$-clean by \cite[Example 2.3]{FY}.
\end{exa}

 However, for some type of polynomials, strong cleanness and strong $g(x)$-cleanness are equivalent.

\begin{thm}\label{thm:12}
Let $R$ be a ring and $g(x) \in (x-a)(x-b)C(R)[x]$ with $a, b\in
C(R)$. Then the following hold:
\begin{enumerate}
\item $R$ is strongly $(x-a)(x-b)$-clean if and only if $R$ is strongly clean and $(b-a) \in U(R)$.
\item If $R$ is strongly clean and $(b-a) \in U(R)$, then $R$ is strongly $g(x)$-clean.
\end{enumerate}
\end{thm}
\begin{proof} $(1).$  ``$\Leftarrow $''. Let $r \in R$. Since $R$ is strongly clean and $(b-a) \in U(R)$, $\frac{r-a}{b-a}=e+u$
 where $e^2=e \in R,$ $u \in U(R)$, and $eu=ue$.  Thus, $r=[e(b-a)+a]+u(b-a)$ where $u(b-a) \in U(R), [e(b-a)+a-a][e(b-a)+a-b]=0$,
 and $[e(b-a)+a]u(b-a)=u(b-a)[e(b-a)+a]$.
 Hence, $R$ is strongly $(x-a)(x-b)$-clean.
\medskip

``$\Rightarrow $''. Since $a$ is strongly $(x-a)(x-b)$-clean, there
exist $u \in U(R)$ and $s \in R$ such that $a=s+u$ with
$(s-a)(s-b)=0$ and $su=us$. Hence, $s=b$. So $(b-a)\in U(R)$. Let $r
\in R$. Since $R$ is strongly $(x-a)(x-b)$-clean, $r(b-a)+a=s+u$
where $(s-a)(s-b)=0,$ $u \in U(R)$, and $su=us$. Thus,
$r=\frac{s-a}{b-a}+\frac{u}{b-a}$ where $\frac{u}{b-a} \in U(R)$,
$(\frac{s-a}{b-a})^2=\frac{(s-a)(s-b+b-a)}{(b-a)^2}=
\frac{(s-a)(b-a)}{(b-a)^2}=\frac{s-a}{b-a}$, and
$\frac{s-a}{b-a}\cdot\frac{u}{b-a}=\frac{u}{b-a}\cdot\frac{s-a}{b-a}$.
So $R$ is strongly clean.
\medskip

$(2).$ By $(1)$.
\end{proof}
\begin{cor}\label{cor:13}
For a ring  $R$, $R$ is strongly clean if and only if $R$ is
strongly $(x^{2}+x)$-clean.
\end{cor}
\begin{proof}
 It follows from Theorem \ref{thm:12} by letting $a=0$ and $b=-1$.
\end{proof}
\begin{rem}\label{rem:14}
The equivalence of strong $(x^{2}+x)$-cleanness and strong cleanness
is a ring property since it holds for a ring $R$ but it may fail for
a single element. For example, $1+1=2\in {\mathbb Z}$ is
 strongly clean but $2$ is not strongly $(x^{2}+x)$-clean in ${\mathbb
 Z}$.
\end{rem}
\begin{exa}\label{exa:2.7} Let $C(X)$ denote the ring of all real-valued continuous functions
from a topological space $X$ to the real number field ${\mathbb R}$
 and $C^{*}(X)$ denote the subring of $C(X)$ consisting of
all bounded functions in $C(X)$. If $X$ is strongly
zero-dimensional, then $C(X)$ and $C^{*}(X)$ are strongly
$(x^2-nx)$-clean rings for any $n\in {\mathbb N}$ since $C(X)$ and
$C^{*}(X)$ are strongly clean \cite{Az02, Mc03} and $n$ is
invertible in $C(X)$ and $C^{*}(X)$. If $X$ is a $P$-space, then
${\mathbb M}_{k}(C(X))$ is strongly $(x^2-nx)$-clean for any $n, k
\in {\mathbb N}$ because ${\mathbb M}_{k}(C(X))$ is strongly clean
by \cite{FY2}.
\end{exa}
\section{general properties of strongly $g(x)$-clean rings}
Let $R$ and $S$ be rings and $\theta: C(R)\rightarrow C(S)$ be a
ring homomorphism with $\theta (1)=1$. For $g(x)=\Sigma a_ix^i\in
C(R)[x]$, let $\theta^{'}\left(g(x)\right):=\Sigma \theta(a_i)x^i\in
C(S)[x]$. Then $\theta$ induces a map $\theta^{'}$ from $C(R)[x]$ to
$C(S)[x]$. If $g(x)$ is a polynomial with coefficients in ${\mathbb
Z}$, then $\theta^{'} (g(x))=g(x)$.
\begin{prop}\label{prop:4} Let $\theta: R \rightarrow S$ be a ring epimorphism. If
$R$ is strongly $g(x)$-clean, then $S$ is strongly
$\theta^{'}(g(x))$-clean.
\end{prop}
\begin{proof}
Let $g(x)=a_{0}+a_{1}x+\cdots+a_{n}x^{n}\in C(R)[x]$. Then
$\theta^{'}
(g(x))=\theta(a_{0})+\theta(a_{1})x+\cdots+\theta(a_{n})x^{n}\in
C(S)[x]$. For any $s \in S$, there exists $r \in R$ such that
$\theta(r)=s$. Since $R$ is strongly $g(x)$-clean, there exist $t
\in R$ and $u \in U(R)$ such that $r=t+u$ with $g(t)=0$ and $tu=ut$.
Then $s=\theta(r)=\theta(t)+\theta(u)$ with $\theta(u)\in U(S)$,
$\theta^{'}(g(x))|_{x=\theta(t)}=0$, and
$\theta(t)\theta(u)=\theta(u)\theta(t)$. So $S$ is strongly
$\theta^{'}(g(x))$-clean.
\end{proof}
\begin{cor} If $R$ is $g(x)$-clean, then for any ideal $I$ of $R$,
$R/I$ is $\overline{g}(x)$-clean with $\overline{g}(x) \in
C(R/I)[x]$.
\end{cor}
\begin{cor}\label{cor:7.1}
Let $R$ be a  ring, $g(x)\in C(R)[x]$, and $1<n \in {\mathbb N}$. If
${\mathbb T}_{n}(R)$ is strongly $g(x)$-clean, then $R$ is strongly
$g(x)$-clean.
\end{cor}
\begin{proof}
Let $A=(a_{ij}) \in {\mathbb T}_{n}(R)$ with $a_{ij}=0$ and $1\leq
j<i \leq n$. Note that $\theta: {\mathbb T}_{n}(R) \rightarrow R$
with $\theta(A)=a_{11}$ is a ring epimorphism.
\end{proof}
\begin{cor}\label{cor:7}
Let $R$ be a ring and $g(x)\in C(R)[x]$. If the formal power series
ring $R[[t]]$ is strongly $g(x)$-clean, then  $R$ is strongly
$g(x)$-clean.
\end{cor}
\begin{proof} This is because $\theta: R[[t]]\rightarrow R$ with $\theta(f)=a_{0}$ is a ring
epimorphism where $f=\sum_{i\geq 0} a_{i}t^{i}\in R[[t]]$.
\end{proof}
\begin{prop}\label{prop:6}
Let $g(x)\in {\mathbb Z}[x]$ and $\{R_{i}\}_{i\in I}$ be a family of
rings. Then $\prod_{i \in I}R_{i}$ is strongly $g(x)$-clean if and
only if $R_{i}$ is strongly $g(x)$-clean for each $i \in I$.
\end{prop}
\begin{proof} It is clear by the definition and Proposition
\ref{prop:4}.
\end{proof}

For strongly clean rings, the author in \cite{CYZ051,C06,C02} proved
that if $R$ is a strongly clean ring and $e^{2}=e \in R$, then the
corner ring $eRe$ is strongly clean. For strongly $g(x)$-clean
rings, we have the following result:

\begin{thm} Let $R$ be a strongly $(x-a)(x-b)$-clean ring with
$a, b\in C(R)$. Then for any $e^{2}=e\in R, eRe$ is strongly
$(x-ea)(x-eb)$-clean. In particular, if $g(x) \in
(x-ea)(x-eb)C(R)[x]$ and $R$ is strongly $(x-a)(x-b)$-clean with $a,
b \in C(R)$, then $eRe$ is strongly $g(x)$ -clean.
\end{thm}
\begin{proof} By Theorem \ref{thm:12}, $R$ is
strongly $(x-a)(x-b)$-clean if and only if $R$ is strongly clean and
$b-a \in U(R)$. If $R$ is strongly clean, then $eRe$ is strongly
clean by \cite{C02}. Again by Theorem \ref{thm:12}, $eRe$ is
strongly $(x-ea)(x-eb)$-clean.
\end{proof}

However, generally, strongly $g(x)$-clean property is not a Morita
invariant: When $g(x)=(x-a)(x-b)$ where $a, b \in C(R)$ with $b-a
\in U(R)$, the matrix ring over the local ring ${\mathbb Z}_{(p)}$
is not strongly clean \cite{CYZ051} (hence, not strongly
$g(x)$-clean).

\section{strongly $g(x)$-clean rings vs rings generated by units and roots of $1$}
For any $n \in {\mathbb N}$, $U_{n}(R)$ denotes the set of elements
of $R$ which can be written as a sum of no more than $n$ units of
$R$ \cite{H74}. A ring $R$ is called \emph{generated by its units}
if $R=\bigcup_{n=1}^{\infty}U_{n}(R)$. We use strong
$g(x)$-cleanness to characterize some rings in which every element
can be written as the sum of unit and a root of $1$ which commute.
\begin{thm}\label{thm:16} Let $R$ be a ring
and $n\in {\mathbb N}$. Then the following are equivalent:
\begin{enumerate}
\item $R$ is strongly $(x^2-2^{n}x)$-clean.
\item $R$ is strongly $(x^2-1)$-clean.
\item $R$ is  strongly clean and $2\in U(R)$.
\item $R=U_{2}(R)$ and for any $a\in R$, $a$ can be expressed as $a=u+v$ with some $u,v\in U(R), uv=vu$, and
$v^{2}=1$.
\end{enumerate}
\end{thm}
\begin{proof}
$(1)\Rightarrow (3)$. To prove $2 \in U(R)$. Suppose $2 \notin
U(R)$, then $\overline{R}=R/(2^{n}R)\neq 0$. Let $2^{n}=s+u$ with
$s^{2}-2^{n}s=0, u \in U(R)$, and $su=us$.
$\overline{0}=\overline{2^{n}}=\overline{s}+\overline{u}$ implies
that $\overline{s}=-\overline{u} \in U(\overline{R})$. But
$\overline{s}^{2}=\overline{s^{2}}=\overline{2^{n}s}=\overline{0}$,
a contradiction. So $2 \in U(R)$. Let $a=0$ and $b=2^{n}$. Then by
$(1)$ of Theorem \ref{thm:12}, $R$ is strongly clean.
\medskip

$(3)\Rightarrow(1)$. By $(1)$ of Theorem \ref{thm:12}, $R$ is
strongly $(x^2-2^{n}x)$-clean.
\medskip

$(3)\Rightarrow (4)$. Let $a \in R$. By ``$(1)\Leftrightarrow(3)$'',
let $n=1$. Then $1-a=s+u$ where $s^2=2s, u \in U(R)$, and $su=us$.
Then $a=(-u)+(1-s)$ with $-u \in U(R), (1-s)^2=1$, and
$(-u)(1-s)=(1-s)(-u)$.
\medskip

$(4)\Rightarrow (3)$. Let $a \in R$. By $(4)$, $1-a=u+v$ where $u
\in U(R), v^2=1$, and $uv=vu$. Thus, $a=(-u)+(1-v)$ with $-u\in
U(R), (1-v)^2=2(1-v)$, and $(-u)(1-v)=(1-v)(-u)$. By
``$(1)\Leftrightarrow(3)$'' and $n=1$, we proved that $(4)$ implies
$(3)$.
\medskip

$(2)\Rightarrow(4)$. If $R$ is strongly $(x^2-1)$-clean, then for
any $r\in R$, there exist $v,u \in U(R)$ such that $r=v+u$ with
$v^2=1$ and $uv=vu$.
\medskip

 $(4)\Rightarrow(2)$. Let $a \in R$. Then $a$ can be expressed as $a=u+v$ with $u,v \in U(R), v^{2}=1$, and $uv=vu$. So $v$ is the root of $x^2-1$.
 Hence, $R$ is
 strongly $(x^2-1)$-clean.
\medskip
\end{proof}
\begin{exa} Rings in Example \ref{exa:2.7} are strongly
$(x^2-nx)$-clean. In particular, they are strongly
$(x^2-2^{n}x)$-clean rings in which every element can be written as
the sum of a unit and a square root of $1$ which commute.
\end{exa}
\begin{exa} Let $F$ be a field and $V$ be a vector space over $F$ of
infinite dimension, and let $R$ be the subring of $E=End_{F}(V)$
generated by the identity and the finite rank transformations. Then
$R$ is strongly clean \cite[Example 7]{C02}. In fact, $E$ is locally
Artinian. So the matrix ring ${\mathbb M }_{k}(E)$ is strongly
$(x^2-nx)$-clean with {\rm char}$F$$\not|n$. If {\rm char}$F\not=2$,
then every element in the matrix ring can be written as the sum of a
unit and a square root of $1$ which commute.
\end{exa}
\begin{exa} Let $A=F[x_1, x_2, \cdots]$ be the polynomial ring in a
countably infinite set of indeterminates $(x_1,  x_2, \cdots)$ over
a field $F$, and let $I=(x_1^{k_1}, x_{2}^{k_2}, x_{3}^{k_3}, ...)$
with $k_{i}>0$. Then $R=A/I$ is a local ring of dimension $0$ which
is not Noetherian. But $R$ is locally Artinian. So the matrix ring
${\mathbb M }_{k}(R)$ is strongly $(x^2-nx)$-clean with {\rm
char}$F$$\not|n$. If {\rm char}$F\not=2$, then every element in the
matrix ring can be written as the sum of a unit and a square root of
$1$ which commute.
\end{exa}
\begin{prop}\label{prop:20}
Let $R$ be a ring with $c, d \in C(R)$ and $d \in U(R)$. If $R$ is
strongly $(x^{2}+cx+d)$-clean, then $R=U_{2}(R)$. In particular, if
$R$ is strongly $(x^{2}+x+1)$-clean, then $R=U_{2}(R)$ is strongly
$(x^{4}-x)$-clean with every element is the sum of a unit and a
cubic root of $1$ which commute with each other.
\end{prop}
\begin{proof}
The first statement is clear. Let $r \in R$. Then $r=s+u$ with $u
\in U(R)$, $s^{2}+s+1=0$, and $su=us$. So $s^{4}-s=0$. Thus, $R$ is
strongly $(x^{4}-x)$-clean. Moreover, ever element in strongly
$(x^{2}+x+1)$-clean ring $R$ can be written as the sum of a unit and
a cubic root of 1 which commute with each other.
\end{proof}
\begin{lem}\cite{FY}\label{lem:11} Let $a\in R$. The following are equivalent for $n \in \mathbb{N}$:
\begin{enumerate}
\item $a=a(ua)^n$ for some $u\in U(R)$.
\item $a=ve$ for some $e^{n+1}=e$ and some $v\in U(R)$.
\item $a=fw$ for some $f^{n+1}=f$ and some $w\in U(R)$.
\end{enumerate}
\end{lem}

\begin{prop}
Let $R$ be an strongly $(x^n-x)$-clean ring where $n \ge 2$ and
$a\in R$. Then either (i) $a=u+v$ where $u\in U(R), v^{n-1}=1$, and
$uv=vu$ or (ii) both $aR$ and $Ra$ contain non-trivial idempotents.
\end{prop}
\begin{proof}  Since $R$ is strongly $(x^n-x)$-clean, $a=s+u$ with $u \in U(R), s^{n}=s$, and $su=us$.
Then $s^{n-1}a=s^{n-1}u+s$. So $(1-s^{n-1})a=(1-s^{n-1})u$. Since
$1-s^{n-1}$ is an idempotent, by Lemma \ref{lem:11},
$(1-s^{n-1})u=vg$ where $v \in U(R)$ and $g^2=g \in R$. So
$g=v^{-1}(1-s^{n-1})a \in Ra$. Suppose (i) does not hold, then
$1-s^{n-1}\not= 0$, this implies $g \not= 0$. Thus, $Ra$ contains a
non-trivial idempotent. Similarly, $aR$ contains a non-trivial
idempotent.
\end{proof}
Finally, we give a property which has nothing to do with rings
generated by units but it relates to strongly $(x^n-x)$-clean rings.
\begin{prop}\label{prop:22}
Let $R$ be a ring and $n\in {\mathbb N}$. Then $R$ is strongly
$(ax^{2n}-bx)$-clean if and only if $R$ is strongly
$(ax^{2n}+bx)$-clean.
\end{prop}
\begin{proof}
$``\Rightarrow".$ Suppose $R$ is strongly $(ax^{2n}-bx)$-clean. Then
for any $r \in R$, $-r=s+u$, $as^{2n}-bs=0, u \in U(R)$, and
$su=us$. So $r=(-s)+(-u)$ where $(-u) \in U(R), a(-s)^{2n}+b(-s)=0$,
and $(-s)(-u)=(-u)(-s)$. Hence, $r$ is strongly
$(ax^{2n}+bx)$-clean. Therefore, $R$ is strongly
$(ax^{2n}+bx)$-clean.
\medskip

$``\Leftarrow".$ Suppose $R$ is strongly $(ax^{2n}+bx)$-clean. Let
$r \in R$. Then there exist $s$ and $u$ such that $-r=s+u$,
$as^{2n}+bs=0, u \in U(R)$, and $su=us$. So $r=(-s)+(-u)$ satisfies
$a(-s)^{2n}-b(-s)=0, -u \in U(R)$, and $(-s)(-u)=(-u)(-s)$. Hence,
$R$ is strongly $(ax^{2n}-bx)$-clean.
\end{proof}
For $2n+1 \in {\mathbb N}$, we do not know if the strong
$(x^{2n+1}-x)$-cleanness of $R$ is equivalent to the strong
$(x^{2n+1}+x)$-cleanness of $R$.
\section*{Acknowledgments}
The first author is partially supported by NSERC, Canada and the
second is supported by the Research Grant from Harbin Institute of
Technology.

\end{document}